\documentclass[12pt]{article} 

\usepackage[utf8]{inputenc} 
\hoffset =- .6cm
\usepackage{graphicx} 

\usepackage{booktabs} 
\usepackage{array} 
\usepackage{paralist}
\usepackage{verbatim}
\usepackage{subfig}
\usepackage{amsmath}
\usepackage{amssymb}

\def\Log{{\rm Log}}

\def\epsilon{\varepsilon}
\def\Log{{\mathrm{Log}}}
\def\e{\varepsilon}

\let\e=\varepsilon

\let\.=\cdot
\let\0=\emptyset

\def\di{\displaystyle}

\def\RR{\mathbb R}

\let\di=\displaystyle
\def\trait (#1) (#2) (#3){\vrule width #1pt height #2pt depth #3pt}

\newcommand{\SE}{\setcounter{equation}{0} \section}
\newcommand{\beq}{\begin{equation}}
\newcommand{\eeq}{\end{equation}}
\newcommand{\baa}{\begin{array}}
\newcommand{\eaa}{\end{array}}
\newcommand{\ba}{\begin{eqnarray}}
\newcommand{\ea}{\end{eqnarray}}

\newtheorem{theorem}{Theorem}[section]

\newtheorem{proposition}[theorem]{Proposition}
\newtheorem{lem}[theorem]{Lemma}
\newtheorem{lemma}[theorem]{Lemma}

\title{\bf Gradient estimates and symmetrization for Fisher-KPP front
propagation with fractional diffusion}
 \author{Jean-Michel {\sc Roquejoffre}$^{\hbox{a }}$, Andrei {\sc
Tarfulea}$^{\hbox{b }}$\\
 \footnotesize{$^{\hbox{a }}$ Institut de Math\'ematiques de Toulouse,
 Universit\'e Paul Sabatier}\\
 \footnotesize{118 route de Narbonne, F-31062 Toulouse Cedex 4, France}\\
 \footnotesize{$^{\hbox{b }}$
Department of Mathematics, Princeton University}\\
 \footnotesize{Fine Hall, Washington Road
Princeton, NJ 08544-1000}\\
 }

\date{}

\begin{document}

\maketitle
\begin{abstract}
In this paper, we study gradient decay estimates for solutions to the multi-dimensional Fisher-KPP equation
with fractional diffusion. It is known that this equation exhibits exponentially advancing
level sets with strong qualitative upper and lower bounds on the solution. 
However, little has been shown concerning the gradient of the solution. 
We prove that, under mild conditions on the initial data, the first and second
derivatives of the solution obey a comparative exponential decay in time.  
We then use this estimate to prove a symmetrization result, which shows
that the reaction front flattens and quantifiably circularizes, losing
its initial structure. 

\end{abstract}

\SE{Introduction}
The goal of this paper is to understand a strong symmetrisation phenomenon,
observed in \cite{C}, for the level sets of reaction-diffusion equations of 
Fisher-KPP type with fractional diffusion. The model under consideration is
\begin{equation}
\partial_t u + \Lambda^\alpha u = f(u) \label{FKPP}
\end{equation}
on ${\mathbb R}^d$. Here, $\alpha\in(0,2)$ and $\Lambda^\alpha$ is the
fractional Laplacian of order $\alpha/2$:
$$(-\Delta)^{\alpha/2}u(x) ={\mathcal F}^{-1}(\vert\xi\vert^{\alpha}\hat u
(\xi))(x)= \displaystyle\lim_{\epsilon \to{}0}{c_{N,\alpha}
\int_{\left|y\right|>\epsilon}\frac{u(x)-u(x+y)}
{\left|y\right|^{d+\alpha}}dy}.
$$
The nonlinear term $f$ is assumed to be of KPP-type \cite{KPP}: $f(0)=f(1)=0$,
$f(y)>0$ for $y \in (0,1)$, and $f'(y) < f(y)/y$. Notation-wise, let $\kappa =
f'(0)$ and $\lambda = \kappa / (d+\alpha)$. 

\noindent Under these assumptions, and for a compactly supported initial datum
$u(.,0):=u_0$, it is well-known that the level sets of $u$ will spread
exponentially fast in time
(Cabré-Roquejoffre \cite{CR}):
\begin{theorem} 
Under the above assumptions,
\begin{itemize}
\item For all $c>\lambda$, we have $\di\lim_{t\to+\infty}\sup_{\vert
x\vert\geq
e^{c t}}u(x,t)=0.$\\
\item For all $c<\lambda$, we have 
$\di\lim_{t\to+\infty}\inf_{\vert x\vert\leq e^{c t}}u(x,t)=1$.
\end{itemize}
\end{theorem}
Thus, when renormalized by the exponential, the level sets of $u$ asymptotically
look round. It is then natural to ask whether this property holds in a more
precise fashion, and a first answer is that given in Cabré-Coulon-Roquejoffre
\cite{CCR}: for a given $h\in(0,1)$ we have
$$
\{u=h\} \subseteq \{C^{-1}e^{\lambda t}\leq\vert x\vert\leq Ce^{\lambda t}\},
$$
for some constant $C>0$. The next question is whether one can make this constant
more precise. Given the rapid growth of the level sets, one could expect the
possibility for very erratic behavior. To confirm this, the case was
investigated numerically by A.-C. Coulon in her PhD thesis. We reproduce here a
sample of her simulations. The initial datum (Fig. 1) is pictured below,
\begin{figure}[h]
\centering
		\includegraphics[width=6cm,height=4cm]{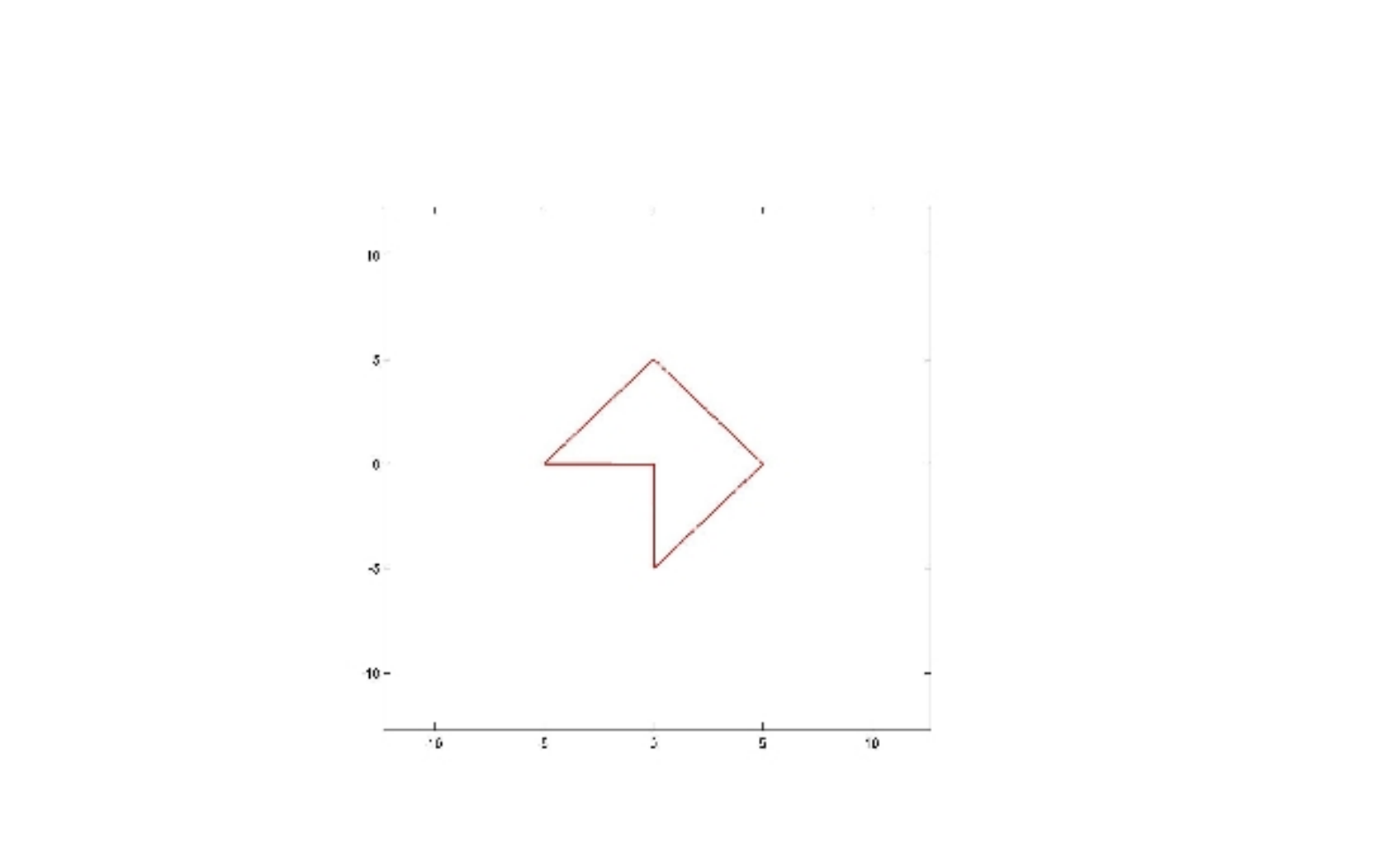}	
		\caption{Level sets of $u_0$}			
\end{figure}

\noindent and the time evolution is shown in Fig. 2.
\begin{figure}[h]
	\centering
		\includegraphics[width=8cm,height=4cm]{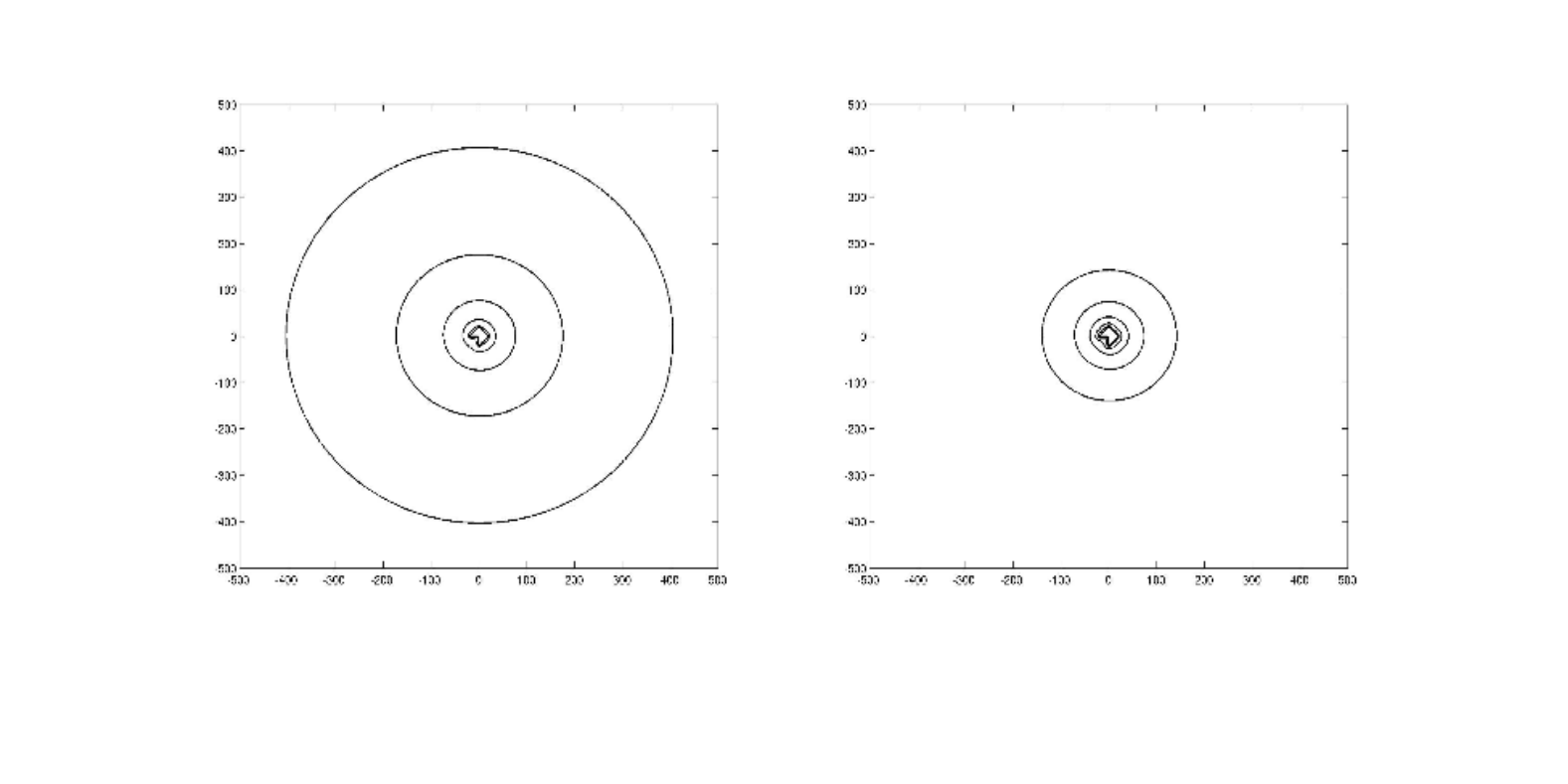}	
		\caption{Level sets of $u$: $\alpha=1$ (left) and $\alpha=1.6$
(right)}		
\end{figure}
\noindent  Surprisingly, this strongly suggests asymptotic symmetrization.
So we may ask whether a very general theorem of Jones \cite{CKRTJ} applies. It
concerns  the solutions of 
(\ref{FKPP}) with $\alpha=2$:
\begin{equation}
\partial_t u -\Delta u = f(u) \label{e1.6}
\end{equation}
with any nonlinearity $f$ - not limited to KPP.  The theorem asserts that, if
$u_0$ is compactly supported, then, for any later time, and for any regular
value $h$ of $u$, and for
any $x$ on the level set $\{u(.,t)=h\}$, the normal line $\Delta_{x,t}$ to the
level set passing through $x$ intersects the convex hull of
${\mathrm{supp}}\, u_0$. If the level sets of $u$ expand, this has very important
consequences: the level sets of $u$ symmetrize asymptotically, and have bounded
oscillation. This is a very remarkable result, one reason being that it does not
depend on the precise expansion rate of the level sets of $u$. In particular, if
the nonlinearity $f$ is of KPP type, then
\begin{equation}
\label{e1.7}
\{u=h\}=\{\vert x\vert=2\sqrt\kappa t-\frac{d+2}2{\Log t}+O(1)\}.
\end{equation}
The term $2\sqrt\kappa t$ is due to Aronson-Weinberger \cite{AW}, and the
logarithmic correction to Bramson \cite{Bram} ($d=1$), and G\"artner \cite{G}
($d>1$).\\

\noindent  Let us briefly describe a proof, due to H. Berestycki \cite{B},  of
Jones's theorem: assume the contrary, then there is a hyperplane $H$ containing
$\Delta_{x,t}$ and such that the convex hull of ${\mathrm{supp}}\, u_0$ lies
strictly in $H^-$, the lower half space bounded by $H$. If $\hat u$ is the
reflection of $u(.,t)$ about $H$, and $v=u-\hat u$, then $v$ satisfies a linear
equation and is positive in $H^-$, since it is positive at $t=0$ and vanishes
on $H$. At time $t$, any derivative of $v$ in a direction normal to $H$ is
nonzero, contradicting that $\Delta_{x,t}\subset H$. We immediately see that
this argument cannot be applied to our case, simply because $v$ would need to be
nonnegative in $H^+$ instead of $H^-$,  which is impossible.\\

\noindent We are going to show that, nevertheless, a rather strong form of
symmetrization occurs. The ingredient  is the following gradient estimate for
$u$; we believe that it is of general interest.
\begin{theorem}\label{t1.1}
Assume $u(x,0):=u_0(x)$ to be continuous, nonnegative and nonzero, and
$u_0(x)=O(e^{-\e\vert x\vert})$ as $\vert x\vert\to+\infty$, for some $\e>0$. 
Then we have, for a solution $u(x,t)$ of \eqref{FKPP}, a
universal constant $C$, and a $\delta > 0$ (depending on $\lambda$ and $\alpha$):
\begin{equation}
\label{e1.1}
\vert \nabla u(x,t)\vert\leq C e^{-\delta t}u(x,t),
\end{equation}
and
\begin{equation}
\label{e1.2}
\vert \nabla^2 u(x,t) \vert \leq C e^{-\delta t}u(x,t).
\end{equation}
\end{theorem}
Theorem \ref{t1.1} allows us to prove an analogous estimate for the fractional
Laplacian. This permits us in turn to reduce the problem to the ODE $\dot
u=(1+O(e^{-\delta t}))f(u)$, which is much simpler to study. The main result of
the paper is
\begin{theorem}\label{t1.2}
Let $u_0$ be as in Theorem \ref{t1.1}. For every $h\in(0,1)$, there is a constant $q_h
(u_0)>0$ and $\delta>0$ such that
\begin{equation}
\label{e1.5}
\{u=h\}=\{\vert x\vert=q_h (u_0)e^{\lambda t}(1+O(e^{-\delta t}))\}.
\end{equation}
\end{theorem}
This provides an explanation to the observed behavior of $u$. However we note
that our results could be improved in two ways.
\begin{itemize}
\item  The assumptions on $u_0$ seem slightly non optimal: indeed, it would
certainly suffice to have $u_0(x)=O(\vert x\vert^{d+\alpha+\epsilon})$ with
$\epsilon>0$. Our theorem would remain valid under that assumption, at the cost
of heavier computations. However, when $u_0$ decays like (or slower than) $\vert
x\vert^{d+\alpha}$ different phenomena occur, as was observed in \cite{C}.  
\item We do not go as far as proving a Jones type theorem. Indeed what we obtain
is a precise estimate of the normalized level sets instead of the true level sets. In other
words, lower order (but still exponentially growing) terms may prevent a
full symmetrization.
\end{itemize}
\noindent These issues  will be investigated in a future paper.\\

\noindent In Section 2, we gather some known (but useful)
facts that will be used throughout the proofs of Theorems \ref{t1.1} and
\ref{t1.2}. Section 3 is devoted to the proof of Theorem \ref{t1.1}, while Section 4
is devoted to estimating the fractional Laplacian. Theorem \ref{t1.2} is then proved
in Section 5.

\SE{Preliminary material}
 The gradient estimate \eqref{e1.1} will be obtained 
by examining its representation in three different ranges, which is reflected by
the collection of results below, that we recall for the reader's convenience.
The proof for the estimate on second derivatives \eqref{e1.2} follows a similar
approach, but will itself make use of \eqref{e1.1}.
\subsection{Invariant coordinates}
The starting point is the following estimate,  proved in \cite{CCR}
\begin{theorem}
\label{t2.1}
We have, for a universal $C>0$:
\begin{equation}
\frac{C^{-1}}{1 + e^{-\kappa t} \left| x \right|^{d+\alpha}} \leq u(x,t) \leq
\frac{C}{1 + e^{-\kappa t} \left| x \right|^{d+\alpha}}. \label{u-bounds}
\end{equation}
\end{theorem}
This motivates the introduction of  the invariant coordinates $\xi = x
e^{-\lambda t}$:
\begin{equation}
\partial_t u - \lambda \xi \cdot \nabla_\xi u + e^{-\alpha \lambda t}
\Lambda^\alpha u - f(u) = 0. \label{xi-FKPP}
\end{equation}
\noindent{}For any fixed coordinate $x_i$, set $\phi = \partial_{x_i} u$. We
then have
\begin{equation}
\label{d-FKPP}
\partial_t \phi + \Lambda^\alpha \phi = f'(u) \phi.
\end{equation}
Letting $v = \partial_{\xi_i} u = e^{\lambda t} \phi$ yields
\begin{equation}
\partial_t v - \lambda \xi \cdot \nabla v + e^{-\alpha \lambda t} \Lambda^\alpha
v - f'(u)v - \lambda v = 0. \label{d-xi-FKPP}
\end{equation}
For another fixed coordinate $x_j$, set $\psi = \partial_{x_j} \phi =
\partial_{x_j} \partial_{x_i} u$ and $\bar{\phi} = \partial_{x_j} u$, with $\bar{v} = e^{\lambda t}
\bar{\phi}$. This then
gives us
\begin{equation}
\partial_t \psi + \Lambda^\alpha \psi = f''(u) \phi \bar{\phi} + f'(u) \psi ,
\label{dd-FKPP}
\end{equation}
along with the associated equation in exponential coordinates: for $V =
\partial_{\xi_j} \partial_{\xi_i} u$, we have
\begin{equation}
\label{dd-xi-FKPP}
\partial_t V - \lambda \xi \cdot \nabla V + e^{-\alpha \lambda t}
\Lambda^{\alpha} V - f'(u)V - f''(u) v \bar{v} - 2 \lambda V = 0.
\end{equation}
A consequence of Theorem \ref{t1.2} is the large time convergence of the
initially compactly supported (or initially exponentially decreasing) solutions
of \eqref{xi-FKPP} to a steady profile. Notice indeed that the
steady equation
\begin{equation}
 - \lambda \xi \cdot \nabla_\xi u =f(u) \label{e2.5}
\end{equation}
has a unique one-parameter family of radial solutions $(u_\tau(\vert
\xi\vert))_{\tau>0}$. Then, assuming Theorem \ref{t1.2}:
\begin{theorem}
\label{t2.1}
There exists a $\tau_\infty(u_0)>0$ such that
$$
\lim_{t\to+\infty}u(\xi,t)=u_{\tau_\infty(u_0)}(\xi),
$$
uniformly on compact subsets of $\RR^d$.
\end{theorem}
\subsection{Heat kernel}
Let $\rho_\alpha(x,t)$ be the heat kernel of $\Lambda^\alpha$, in other words
the solution of $\rho_t+\Lambda^\alpha\rho=0$ having, as initial datum, the
Dirac mass at 0.
\begin{proposition}
\label{p2.1}
There is $C>0$ such that, for large $t$ and $x$, we have
$\rho_\alpha(x,t)=t^{-d/\alpha}p_\alpha(x/t^{d/\alpha})$, with
\begin{equation}
\label{e2.1}
 \frac{C^{-1}}{\vert \zeta\vert^{d+\alpha+1}}\leq \vert \nabla
p_\alpha(\zeta)\vert\leq\frac{C}{\vert \zeta\vert^{d+\alpha+1}}, 
\end{equation}
and there is $c_\alpha>0$, $\delta>0$ such that
\begin{equation}
\label{e2.2}
 \frac{C^{-1}}{\vert \zeta\vert^{d+\alpha+\delta}}\leq
p_\alpha(\zeta)-\frac{c_\alpha}{\vert\zeta\vert^{d+\alpha}}\leq\frac{C}{\vert
\zeta\vert^{d+\alpha+\delta}}.
\end{equation}

\end{proposition}
It is important to note that the estimates \eqref{e2.1} and \eqref{e2.2} are
intended for $\zeta$ bounded away from $0$; $\rho_\alpha$ and its derivatives
are bounded functions (in $x$) for any $t>0$. A standard way to prove the above
proposition is to write 
$$\rho_\alpha(x,t)= {\mathcal F}^{-1} \left( e^{-\left| \xi \right|^\alpha
t}\right),
$$
and to evaluate the inverse Fourier transform with the aid of Polya integrals;
see for instance \cite{Ko}.
\subsection{Comparison principles}
We will also require the following  easy extension of the Maximum Principle:

\begin{theorem}[\textbf Selective Comparison Principle]\label{SCP}
Let $\Omega(t)$ be a time-dependent family of compact domains in ${\mathbb R}^d$
(in the $\xi$ variable) with smooth boundaries and continuous time dependence;
that is, $ \lbrace (\xi,t) \ | \ \xi \in \Omega(t) \rbrace $ is an open set in
${\mathbb R}^d \times {\mathbb R}_{+}$. Let $v(\xi, t)$ be the solution to
\eqref{d-xi-FKPP} and let $w(\xi, t)$ be a positive function such that $w >
\vert v \vert$ at $t=0$. If, for all $t>0$ we have $w> \vert v \vert$ on the
closure of the complement of $\Omega(t)$ and
$$ \partial_t w - \lambda \xi \cdot \nabla w + e^{-\alpha \lambda t}
\Lambda^\alpha w - f'(u) w - \lambda w > 0 $$
\noindent{}everywhere inside $\Omega(t)$, then $w(\xi, t) \geq \vert v(\xi, t)
\vert $ on all of ${\mathbb R}^d \times {\mathbb R}_{+}$.
\end{theorem}

\noindent{\sc Proof.} We first show that $w>v$ by looking at the equation for $q
= w-v$ on $\Omega(t)$ which, by assumption, starts positive at $t=0$. Observe
that, if $q$ achieves a global minimum value of $0$ at time $\bar{t}$ and
location $\bar{\xi} \in \Omega(\bar{t})$, then we easily have $\partial_t
q(\bar{\xi}, \bar{t}) > 0$ (all other terms on the left are nonpositive); this,
however, crucially requires that $w>v$ outside $\Omega(\bar{t})$. So, by
continuity in time, $q$ can never become negative inside $\Omega(t)$. To
complete the proof, we need to show that $w>-v$. But this is done in precisely
the same manner as before, now examining the equation for $\bar{q} = w+v$ and
showing that it too remains positive for all time. \hfill $\bullet$
\\

\noindent{}Also note that the above result (with a minor modification) also
holds for a solution $V$ of \eqref{dd-xi-FKPP}. That is, if $W$ is a positive
function such that $W>\vert V \vert$ on the closure of the complement of
$\Omega(t)$ for all $t>0$ and also $W> \vert V \vert$ on $\Omega(0)$, and
$$ \partial_t W - \lambda \xi \cdot \nabla W + e^{-\alpha \lambda t}
\Lambda^\alpha W - f'(u) W - 2 \lambda W > \vert f''(u) v \bar{v} \vert $$
everywhere inside $\Omega(t)$, then $W(\xi, t) > \vert V(\xi,t) \vert$ on all of
${\mathbb R}^d \times {\mathbb R}_{+}$. The proof is exactly the same since the
equations for $W-V$ and $W+V$ have nonnegative inhomogeneities on the right side
of the inequality.\\

\noindent{}Finally let us mention the following version of Kato's inequality and
its well-known consequence for the Fisher-KPP equation.
\begin{proposition}
\label{p2.10}
If $u(x)$ is smooth, then 
$$\Lambda^\alpha\vert u(x)\vert \leq {\mathrm{sgn}}(u)\Lambda^\alpha u(x),
$$
in the distributional sense (and in the classical sense if $\alpha<1$).
\end{proposition}
\noindent{\sc Proof.} Recall the elementary inequality $\vert a\vert-\vert
b\vert\geq{\mathrm{sgn}}(b)(a-b)$. This implies, for all
$(x,h)\in\RR^d\times\RR^d$:
$$
\vert u(x+h)\vert-\vert u(x)\vert\geq {\mathrm{sgn}}(u(x))(u(x+h)-u(x)).
$$
Integrating both sides in the variable $h$ over
$\RR^d\backslash(-\varepsilon,\varepsilon)^d$ and letting $\varepsilon\to0$
yields the result. \hfill $\bullet$ \\
\\
\noindent This implies the following lemma:
\begin{lemma}
\label{l2.5}
If $\phi(x,t)$ is either $u_t$, or $\partial_{x_i}u(x,t)$, ($i\in\{1,...,d\}$),
where $u(x,t)$ is the solution of \eqref{FKPP}, then 
$$\partial_t\vert\phi\vert+\Lambda^\alpha\vert\phi\vert\leq
\kappa\vert\phi\vert.
$$
Additionally, we have
$$ \partial_t \vert \psi \vert + \Lambda^\alpha \vert \psi \vert \leq \vert
f''(u) \phi \bar{\phi} \vert + \kappa \vert \psi \vert
$$
for $\psi(x,t)$ a solution to \eqref{dd-FKPP}.
\end{lemma}
\noindent{\sc Proof.} Multiplying \eqref{d-FKPP} by ${\mathrm{sgn}}(\phi)$
(respectively \eqref{dd-FKPP} by ${\mathrm{sgn}}(\psi)$) and using Proposition
\ref{p2.10} yields the result. \hfill$\bullet$
\SE{The gradient decay estimate}

Throughout this section, all inequalities will be up to a constant independent
of the solutions. We also weaken our requirements on the decay of the initial
data. Specifically, we insist that
$$ 0 \leq u_0(x) \leq \frac{C}{1 + \left| x \right|^{d+\alpha+1}}. $$
With the notations of Section 2, we would like to prove an exponential decay
rate for the gradient \eqref{e1.1} and second derivatives \eqref{e1.2}
of the form\\
\begin{equation}
\vert \phi(x,t) \vert \leq u(x,t) e^{-\delta t} \ \ \ \ \text{and} \ \ \ \ \
\vert \psi(x,t) \vert \leq u(x,t) e^{-\delta t}.  \label{GDE}
\end{equation}

\noindent{}To do this, we will need to examine both representations for the
evolution of the derivative (\eqref{d-FKPP} and \eqref{d-xi-FKPP}). We have to
distinguish three different ranges of $\left| x \right|$: long range
($\vert x\vert>\!> e^{\lambda  t}$), medium range ($\vert x\vert\sim e^{\lambda
t}$), and short range ($\vert x\vert <\!<e^{\lambda t}$). The proof of
\eqref{e1.1} is stand-alone and presented fully. The proof of \eqref{e1.2} is
nearly identical (and we present it in parallel), but assumes that \eqref{e1.1}
holds a priori.

\subsection{Long Range: Kernel Estimates}

\begin{lem}
For any fixed $c>0$ and for any $x$ in the range $\left| x \right| \geq
ce^{\lambda (\alpha+1) t}$, a solution $\phi$ of \eqref{d-FKPP} satisfies
$$ | \phi (x,t) | \leq u(x,t) e^{-\alpha \lambda t}. $$
\end{lem} 
\noindent{}Later, we will set the value of $c$ based on the parameters of the problem.
For now, $c$ only affects the final constant that appears in estimate \eqref{GDE}. As
such, we suppress its notation in the proof of the lemma.
We also draw special attention to the fact that, in the long range, we have an explicit
exponent for the decay rate of $\lambda \alpha$. This particular value will later facilitate a
compatibility condition between the long and medium ranges.\\
\\
{\noindent\sc Proof.} We first prove the estimate for $\phi$. Set
$g(u)=f(u)-\kappa u$. By Duhamel's formula we have
\begin{equation}
\label{e3.1}
\begin{array}{rll}
 \phi(x,t)=&e^{\kappa t}\displaystyle\int\partial_{y_i}
\rho_\alpha(x-y,t)u_0(y)dy\\
 +& \displaystyle\int_0^t\int e^{\kappa(t-s)}\rho_\alpha(x-y,t-s) \partial_{y_i}
g(u(y,s))dyds.
 \end{array}
\end{equation}
Note that   $\|\partial_x \rho_\alpha(.,t) \|_{L^1}$  is not integrable near
$t=0$, for $\alpha<1/2$. Hence the second term in the expression for  $\phi$ 
must be split into two pieces; one in a neighborhood of $s = t$ (where
$e^{\kappa (t-s)}$ is essentially constant) and the rest where we may integrate
by parts. In total, we obtain the following upper bound on $\phi(x,t)$:
$ \vert \phi(x, t)\vert \leq I_0 + I_1 + I_2 $,
where
$$ I_0 = e^{\kappa t} \int_{{\mathbb R}^d} u_0(x-y) \left| \partial_y
\rho_\alpha (y,t) \right| dy, $$
$$ I_1 = \int_0^{t-1} e^{\kappa (t-s)} \int_{{\mathbb R}^d} \left|
f(u(x-y,s))-\kappa u(x-y,s) \right| \left| \partial_y \rho_\alpha (y, t-s)
\right| dy ds, $$
$$ I_2 = \int_{t-1}^t \int_{{\mathbb R}^d} \left| f'(u(x-y,s))
\phi(x-y,s)-\kappa \phi(x-y,s) \right| \rho_\alpha (y, t-s) dy ds. $$

\noindent{}So we proceed by estimating each of the three integrals. We begin
with $I_0$ and apply the naive estimate for $\partial_x \rho_\alpha$
\eqref{e2.1} and a concrete decay rate for the initial data $u_0$.
\begin{equation}
\label {I_0}
I_0 \leq e^{\kappa t} P(t) \int_{{\mathbb R}^d} \frac{1}{1+\left| x-y
\right|^{d+\alpha+1}} \frac{1}{1+\left| y \right|^{d+\alpha+1}} dy,
\end{equation}
where $P(t)$ is a polynomial in $t$. The specific form and degree of $P$ are not
needed since we will be getting a slight exponential decay in time, and that
dominates any residual polynomial growth. To handle the main integral in
(\ref{I_0}), we split into two regions: $\lbrace \left| y \right| \leq \left| x
\right| / 2 \rbrace $ and $\lbrace \left| y \right| > \left| x \right| / 2
\rbrace $. Observe that $\left| x-y \right| \geq \left| x \right| - \left| y
\right|$ which, in the first case, means $\left| x-y \right| \geq \left| x
\right| /2 $. Thus

$$ \int_{\left| y \right| \leq \left| x \right| / 2} \frac{1}{1+\left| x-y
\right|^{d+\alpha+1}} \frac{1}{1+\left| y \right|^{d+\alpha+1}} dy \leq
\frac{1}{\left| x \right|^{d+\alpha+1}} \int_{{\mathbb R}^d} \frac{dy}{1+\left|
y \right|^{d+\alpha+1}} \leq \frac{1}{\left| x \right|^{d+\alpha+1}} $$

\noindent{}for $x$ in the long range. For the second case,

$$ \int_{\left| y \right| > \left| x \right| / 2} \frac{1}{1+\left| x-y
\right|^{d+\alpha+1}} \frac{1}{1+\left| y \right|^{d+\alpha+1}} dy \leq
\frac{1}{\left| x \right|^{d+\alpha+1}} \int_{{\mathbb R}^d} \frac{dy}{1+\left|
x-y \right|^{d+\alpha+1}} \leq \frac{1}{\left| x \right|^{d+\alpha+1}}. $$

\noindent{}Finally, observe that condition (\ref{u-bounds}) and the size of
$\left| x \right|$ imply that $u(x,t) \approx e^{\kappa t} / \left| x
\right|^{d+\alpha}$. Since $\left| x \right| \geq e^{\lambda (\alpha+1) t}$, we
then have

$$ I_0 \leq e^{\kappa t} \left| x \right|^{-(d+\alpha+1)} \approx u \left| x
\right|^{-1} \leq u e^{-\lambda (\alpha+1) t}, $$

\noindent{}which gives us our desired comparative decay for $I_0$.\\
\\
Continuing, we need to estimate the factors involving $f$ and $f'$. Since $f$ is
smooth enough to have a partial Taylor expansion, we have that $f(u) = f(0) + f'(0) u +
O(u^2)$. Therefore $f(u) - \kappa u = O(u^2)$, and this (along with
(\ref{u-bounds})) gives one of the upper bounds we will use for $I_1$.
\begin{equation}
I_1 \leq \int_0^{t-1} e^{\kappa(t-s)} Q(t,s) \int_{{\mathbb R}^d} \frac{1}{1 +
e^{-2 \kappa s} \left| y \right|^{2(d+ \alpha)}} \frac{1}{1+\left| x-y
\right|^{d+\alpha+1}} dy ds. \label{I_1}
\end{equation}

\noindent{}Here $Q(t,s)$ is algebraic in $t$ and $s$. As discussed earlier,
$Q(t,s)$ will not be integrable up to $s = t$, but it is at most polynomial in
$t$ on the interval $0 \leq s \leq t-1$. As with $I_0$, we handle the main
integral by splitting into the same two regions: $\lbrace \left| y \right| \leq
\left| x \right| / 2 \rbrace $ and $\lbrace \left| y \right| > \left| x \right|
/ 2 \rbrace $. In the first case, we get
\begin{equation}
\label{I_1a}
\begin{split}
\int_{\left| y \right| \leq \left| x \right| / 2} & \frac{1}{1 + e^{-2 \kappa s}
\left| y \right|^{2(d+ \alpha)}} \frac{1}{1+\left| x-y \right|^{d+\alpha+1}} dy
\\
& \leq \frac{1}{\left| x \right|^{d+\alpha+1}} \int_{{\mathbb R}^d} \frac{dy}{1
+ e^{-2 \kappa s} \left| y \right|^{2(d+ \alpha)}} \approx \frac{e^{\lambda d
s}}{\left| x \right|^{d+\alpha+1}}.
\end{split}
\end{equation}
\noindent{}The last estimate follows from a change of variables. In the second
case, we get
\begin{equation}
\begin{array}{rll}
&\displaystyle \int_{\left| y \right| > \left| x \right| / 2} \frac{1}{1 + e^{-2
\kappa s} \left| y \right|^{2(d+ \alpha)}} \frac{1}{1+\left| x-y
\right|^{d+\alpha+1}} dy\\
 &\displaystyle \leq \frac{e^{2 \kappa s}}{\left| x \right|^{2(d+\alpha)}}
\int_{{\mathbb R}^d} \frac{dy}{1+\left| x-y \right|^{d+\alpha+1}} \approx
\frac{e^{2 \kappa s}}{\left| x \right|^{2(d+\alpha)}} .
 \end{array}
 \label{I_1b}
\end{equation}

\noindent{}Plugging the estimates from (\ref{I_1a}) and (\ref{I_1b}) back into
(\ref{I_1}) shows that

$$ I_1 \leq \int_0^{t-1} Q(t,s) e^{\kappa (t-s)} \left( \frac{e^{\kappa s
\frac{d}{d+\alpha}}}{\left| x \right|^{d+\alpha+1}} +  \frac{e^{2 \kappa
s}}{\left| x \right|^{2(d+\alpha)}} \right) ds \leq P(t) \left( \frac{e^{\kappa
t}}{\left| x \right|^{d+\alpha+1}} + \frac{e^{2 \kappa t}}{\left| x
\right|^{2(d+\alpha)}} \right). $$

\noindent{}Remembering that $u \approx e^{\kappa t} \left| x
\right|^{-(d+\alpha)}$ and $\left| x \right| \geq e^{\frac{1+\alpha}{d+\alpha}
\kappa t}$, we see that

$$ I_1 \leq u P(t) \left( \frac{1}{\left| x \right|} + \frac{e^{\kappa
t}}{\left| x \right|^{d+\alpha}} \right) \leq u P(t) \left( e^{-\lambda
(1+\alpha) t} + e^{-\alpha \kappa t} \right), $$

\noindent{}establishing the comparative exponential decay rate for $I_1$.
Observe that our estimates here would not work if we only had $\left| x \right|
\geq e^{\lambda t}$.\\
\\
We handle $I_2$ in an analogous manner. Note that $f'(u) = f'(0) + O(u)$, so
that $f'(u) \phi - \kappa \phi = O(u \phi)$. But we also have, from
Lemma \ref{l2.5},
$$
\vert\phi(x,t)\vert\leq e^{\kappa t}\rho_\alpha*_{x_i}
\vert\partial_xu_0\vert\leq Ct\frac{e^{\kappa t}}{\vert x\vert^{d+\alpha}},
$$
and so $\vert\phi(x,t)\vert\leq Ctu(x,t)$, recalling that Theorem \ref{t2.1}
implies $u(x,t)\geq Ce^{\kappa t}\vert x\vert^{-(d+\alpha)}$. This again implies
a $O(tu^2)$ upper bound which reduces the estimate on $I_2$ to

$$ I_2 \leq \int_{t-1}^t \bar{Q}(t,s) \int_{{\mathbb R}^d} \frac{1}{1 + e^{-2
\kappa s} \left| y \right|^{2(d+ \alpha)}} \frac{1}{1+\left| x-y
\right|^{d+\alpha}} dy ds. $$

\noindent{}Here, $\bar{Q}(t,s)$ is again algebraic in $t$ and $s$, but can now
be made integrable on $[t-1, t]$ due to the correct time-homogeneity of
$\rho_\alpha (x, t)$. We estimate the convolution integral exactly as we did for
$I_1$ and find that

$$ I_2 \leq \int_{t-1}^t \bar{Q}(t,s) \left( \frac{e^{\kappa s
\frac{d}{d+\alpha}}}{\left| x \right|^{d+\alpha}} + \frac{e^{2 \kappa s}}{\left|
x \right|^{2(d+\alpha)}} \right) ds \leq P(t) \left( \frac{e^{d \lambda
t}}{\left| x \right|^{d+\alpha}} + \frac{e^{2 \kappa t}}{\left| x
\right|^{2(d+\alpha)}} \right). $$

\noindent{}Continuing as we did for $I_1$, we finally get

$$ I_2 \leq u P(t) \left( e^{d \lambda t - \kappa t} + e^{\kappa t - \kappa (1 +
\alpha) t} \right) = u P(t) \left( e^{-\frac{\alpha \kappa t}{d+\alpha}} +
e^{-\alpha \kappa t} \right). $$

\noindent{}Finally, we examine the exponents to conclude that, for $\left| x \right| \geq
e^{\lambda (1+\alpha) t}$, we have $\left| \phi(x,t) \right| \leq u(x,t)
e^{-\lambda \alpha t}$. Here, again, the method would not have worked if we merely
required $\left| x \right| \geq e^{\lambda t}$.\\

\noindent{}Now we turn our attention to $\psi$. Looking at \eqref{dd-FKPP}, we
employ the Duhamel formula to obtain
\begin{equation}
\label{e3.1-d}
\begin{array}{rll}
 \psi(x,t)=&e^{\kappa t}\displaystyle\int\partial_{y_i} \partial_{y_i}
\rho_\alpha(x-y,t) u_0(y)dy\\
 +& \displaystyle\int_0^t\int e^{\kappa(t-s)}\rho_\alpha(x-y,t-s) \partial_{y_i}
g(\phi(y,s))dyds\\
 +& \displaystyle \int_0^t \int e^{\kappa(t-s)} \rho_\alpha(x-y,t-s) f''(u(y,s))
\phi(y,s)^2 dyds.
 \end{array}
\end{equation}

\noindent{}The comparative decay estimates follow analogously to the argument
for $\phi$. The first term has even better decay than $I_0$ since we placed two
derivatives on $\rho_\alpha$. The second term is again split into two parts
$\bar{I}_1$ and $\bar{I}_2$ (for the time intervals $[0,t-1]$ and $[t-1,t]$).
Since we have already concluded that $\vert \phi \vert \leq u e^{-\lambda \alpha t}$
in this range, the analysis for $\bar{I}_1$ proceeds identically.\\

\noindent{}Now, $f'(u)\psi - \kappa \psi = O(u\psi)$. The analysis of
$\bar{I}_2$ would also be identical to that of $I_2$ if we knew that $\vert \psi
\vert = O(u)$ (up to a time-dependent, not exponentially increasing factor). We
prove this by appealing to Lemma \ref{l2.5}.
$$ \vert \psi (x,t) \vert \leq e^{\kappa t} \rho_\alpha *_x \vert \partial_{x_j}
\partial_{x_i} u_0 \vert + \int_0^t e^{\kappa (t-s)} \rho_\alpha (\cdot, t-s)
*_x \vert f''(u) \phi \bar{\phi} \vert(\cdot, s) ds := J_0 + J_1 $$
$J_0$ is bounded by $Ct e^{\kappa t} / \vert x \vert^{d+\alpha}$ as before.
$J_1$ is estimated by the (now proven) comparative exponential bound for $\phi$
(and $\bar{\phi}$) in the long range:
\begin{equation}
\label{e3.2}
\begin{array}{rll}
J_1 \leq & Ce^{-2 \lambda \alpha t} \int_0^t \int e^{\kappa (t-s)} \rho_\alpha (x-y,
t-s) u(y,s)^2 dy ds \\
\leq & Ce^{-2 \lambda \alpha t} \int_0^t e^{\kappa (t-s)} \bar{Q}(t,s) \int \frac{1}{1
+ e^{-2 \kappa s} \left| y \right|^{2(d+ \alpha)}} \frac{1}{1+\left| x-y
\right|^{d+\alpha}} dy ds \\
\leq & P(t)e^{-2 \lambda \alpha t} \left( e^{\kappa t} \vert x \vert^{-(d+\alpha)}
\int_0^t e^{-\alpha \lambda s} ds + e^{\kappa t} \vert x \vert^{-2(d+\alpha)}
\int_0^t e^{\kappa s} ds \right) \\
\leq & P(t) e^{-2 \lambda \alpha t} \left( u(x,t) + u(x,t)^2 \right).
\end{array}
\end{equation}
The third inequality was obtained through the same estimates as for $I_2$.
Notice how $J_1$ inherits an exponential decay from the more rigorous analysis
for $\phi$; using a similar approach for the second term of \eqref{e3.1} would
have merely given us a comparative boundedness for $\phi$ instead of decay.
Hence $J_1$ is dominated by $J_0$. Therefore $\psi = O(tu)$, and we conclude
that the second term of \eqref{e3.1-d} satisfies a comparative exponential decay
rate in the long range.\\

\noindent{}Lastly, observe that the third term of \eqref{e3.1-d} is in fact
dominated by $J_1$. As seen in \eqref{e3.2}, $J_1$ also satisfies a comparative
exponential decay rate (with exponent $-2 \lambda \alpha$). This completes the
mirrored argument, demonstrating that $\vert \psi(x,t) \vert \leq
u(x,t)e^{-\lambda \alpha t}$ as long as $\vert x \vert \geq e^{\lambda (1+\alpha) t}$. \hfill $\bullet$

\subsection{Short Range: stability of the steady state $u\equiv 1$}

The nature of $f$ and condition (\ref{u-bounds}) ensures that there exists a
fixed $\beta > 0$ such that $\lbrace \left| \xi \right| < \frac{3}{2} \beta
\rbrace \subseteq \lbrace f'(u) < -\delta_0 \rbrace$. That is, we have an
exponentially growing ball wherein $u$ is sufficiently close to $1$ that $f'(u)$
is negative by at least a fixed amount. Let $\Omega_0(t) = \lbrace \left| \xi
\right| < \frac{3}{2} \beta \rbrace$. In this region, we use a simple version of
the maximum principle. Assume the positive maximum value of $\phi$ occurs at
time $\bar{t}$ at an interior point $\bar{x} \in \Omega_0(\bar{t})$ (negative
minimum values are handled similarly). Then

$$ \partial_t \phi (\bar{x},\bar{t}) + \Lambda^{\alpha} \phi (\bar{x},\bar{t}) =
f'(u(\bar{x}, \bar{t})) \phi(\bar{x},\bar{t}). $$

\noindent{} However, we have $\Lambda^{\alpha} \phi (\bar{x},\bar{t}) \geq 0.$
Therefore,

$$ \partial_t \phi (\bar{x},\bar{t})  \leq -\delta_0 \phi(\bar{x},\bar{t}).
$$

\noindent{}We conclude that either $\phi$ is decaying exponentially (with
exponent at least $-\delta_0$) on ${\mathbb R}^d$, \emph{or} the size of $\phi$
on $\Omega_0(t)$ is controlled by the size of $\phi$ outside $\Omega_0(t)$.
Equivalently, $v$ is eventually controlled by a fixed constant times
$e^{(\lambda-\delta_0) t}$ whenever $\left| \xi \right| < \frac{3}{2} \beta$ so
long as $v$ obeys the same bound for (say) $\left| \xi \right| > \beta$, which
will be established in the next two sections. Certainly, the previous argument shows
that this is the case in the long range.\\

\noindent To prove the analogous statement for $\psi$, we will assume that \eqref{e1.1}
holds with some uniform $\delta > 0$ in all ranges for $\phi$ (the proof of
which will, as stated, be stand-alone). Since $u$ is bounded by $1$, we have
that
$$ \partial_t \psi(x,t) + \Lambda^\alpha \psi(x,t) \leq -\delta_0 \psi(x,t) + C
e^{-2 \delta t} $$
for all $x \in \Omega_0(t)$. Assume the positive maximum value of $\psi$ occurs
at time $\bar{t}$ at an interior point $\bar{x} \in \Omega_0(t)$. Then
$$ \partial_t \psi(\bar{x},\bar{t}) \leq -\delta_0 \psi(\bar{x}, \bar{t}) + C
e^{-2 \delta \bar{t}} $$
again establishing that \emph{either} $\psi$ is decaying exponentially on
${\mathbb R}^d$, \emph{or} the size of $\psi$ on $\Omega_0(t)$ is controlled by
the size of $\psi$ outside $\Omega_0(t)$. As before, negative minima are treated
similarly.

\subsection{Medium range: selective comparison principle}

For the intermediate range, we appeal to Theorem (\ref{SCP}). Let $\Omega_1(t) =
\lbrace \beta \leq \left| \xi \right| \leq ce^{\lambda \alpha t} \rbrace$, for
$c>0$ to be determined later. We need to find a $w$ such that $w(\xi,0) \geq
\vert v(\xi,0) \vert$, and $w > \vert v \vert$ outside of
$\Omega_1(t)$ for all $t>0$, and such that
\begin{equation}
\partial_t w - \lambda \xi \cdot \nabla w + e^{-\lambda \alpha t} \Lambda^\alpha
w - \lambda w - f'(u) w \geq 0 \label{super}
\end{equation}
\emph{pointwise on $\Omega_1(t)$}. It is important that $w$ be defined on all of
${\mathbb R}^d$ because $\Lambda^\alpha$ is nonlocal; while it can be defined on
compact sets, it becomes a different operator than the one in
(\ref{d-xi-FKPP}).\\
\\
We must therefore provide a candidate supersolution $w$. For $\xi \in
\Omega_1(t)$ and $0< \nu - d-\alpha \leq 1$, define the family of functions
$\tilde{w}_\nu$ by
\begin{equation}
\tilde{w}_\nu(\xi, t) = \frac{\beta^\nu w_0}{\left| \xi \right|^\nu} + C_\nu
e^{-\lambda \alpha t} \left( \frac{\beta^{\nu-d-\alpha}}{\left| \xi \right|^\nu}
- \frac{1}{\left| \xi \right|^{d+\alpha}} \right), \label{W}
\end{equation}
where $C_\nu, w_0 >0$ will be determined later. Observe that each such
$\tilde{w}_\nu$ is radial. To simplify notation, let $y = \left| \xi \right|$.
For fixed $\nu$ and $t$, , $\tilde{w}_\nu$ is the solution of the following ODE
in $y$ (on $\Omega_1(t)$):
\begin{equation}
\label{odeW}
y \frac{d}{dy} \tilde{w}_\nu(y,t) + \nu \tilde{w}_\nu(y,t) = -\bar{C}_\nu
e^{-\lambda \alpha t} y^{-(d+\alpha)}.
\end{equation}
Note that $C_\nu = \bar{C}_\nu /(\nu-d+\alpha)$. $\bar{C}_\nu$ will be
determined below. Technically, $\bar{C}_\nu$ will also depend on $w_0$, but the
dependence will be linear. More importantly, this dependence is not circular:
$w_0$ is the initial condition at $y=\beta$, and this may be fixed before
solving the ODE. $w_0$ itself will depend on the size of $v(\xi, 0)$.\\
\\
Our supersolution must be positive, yet $\tilde{w}_\nu$ is eventually negative
for sufficiently large $y$. We now choose $c$ such that $\tilde{w}_\nu > 0$ on
$\Omega_1(t)$; since $\nu \leq d+\alpha+1$, such a choice is always possible.
Outside of $\Omega_1(t)$, we extend $\tilde{w}_\nu$ to be smooth, bounded,
positive, radial, decreasing in $y$, and $O(y^{-(d+\alpha)})$. This can
obviously be done such that the $L^1$, $C^1$, and $C^2$ norms of $\tilde w$ are
uniformly bounded for all $t>0$. They do, however, depend on $\nu$ and $w_0$. \\
\\
We now define our candidate supersolution as $w(\xi,t) = e^{\eta t}
\tilde{w}_\nu (\left| \xi \right|,t)$ with $\eta$ and $\nu$ fixed; we will
determine their values shortly. Observe that $\xi \cdot \nabla w = y e^{\eta t}
\partial_y \tilde{w}$. To prove (\ref{super}), it therefore suffices to show
that

$$ -e^{\eta t} \partial_t \tilde{w}_\nu - \eta w + \lambda y \partial_y w(y,t) +
\left( \lambda + f'(u(\xi,t)) \right) w(y,t) \leq - e^{-\lambda \alpha t} \left|
\Lambda^\alpha w (y,t) \right| $$

\noindent{}on $\Omega_1$(t). Examining (\ref{W}), we see that $\partial_t
\tilde{w}_\nu > 0$ on our domain, so we get a stronger inequality if we ignore
the first term above; proving the stronger inequality will be sufficient. Since
$f'(u) \leq \kappa = \lambda (d+\alpha)$ and $w > 0$, we can eliminate the $\xi$
dependence entirely and satisfy the above inequality provided that

$$ \lambda y \partial_y w(y,t) + \lambda (d+\alpha+1 -
\frac{\eta}{\lambda})w(y,t) \leq -e^{-\lambda \alpha t} \left| \Lambda^\alpha w
(y,t) \right| .$$

\noindent{}Treating the diffusion as a perturbation term, we need to bound the
size of $\left| \Lambda^\alpha w(y,t) \right|$. Our choice of extension for
$\tilde{w}_\nu$ onto ${\mathbb R}^d$ implies (see \cite{BV}, and estimate (8) of
\cite{CCR}):

$$ \left| \Lambda^\alpha w(y,t) \right| = e^{\eta t} \left| \Lambda^\alpha
\tilde{w}_\nu (y,t) \right| \leq e^{\eta t} \frac{\bar{C}_\nu}{y^{d+\alpha}} $$

\noindent{}on $\Omega_1$(t). Our assumptions on the uniform-in-time bounds of
the norms of $\tilde{w}_\nu$ ensure that the constant $\bar{C}_\nu$ can be made
time-independent (recall that $\nu$ is fixed). It is also linear in $w_0$, as
mentioned before. Thus, proving (\ref{super}) reduces to showing that

$$ y \partial_y w(y,t) + (d+\alpha+1-\frac{\eta}{\lambda})w(y,t) = -e^{\eta
t-\lambda \alpha t} \frac{\bar{C}_\nu}{y^{d+\alpha}}, $$

\noindent{}with $\bar{C}_\nu>0$ now determined. But this is precisely
(\ref{odeW}), multiplied by $e^{\eta t}$ with $\nu = d + \alpha + 1 -
\frac{\eta}{\lambda}$. Hence, (\ref{super}) holds on $\Omega_1$(t). Clearly, we
need $0 \leq \eta < \lambda$. We then choose $w_0$ sufficiently large that $w$
(that is, the function extended to all of ${\mathbb R}^d$) dominates $\vert v
\vert$ at $t = 0$.\\
\\
In order to use Theorem (\ref{SCP}), we need to show that $w- \vert v \vert >0$
outside of $\Omega_1(t)$. It is here that we needed our free parameters $\nu$
and $\eta$; naively, we would like to take $\eta = 0$ and $\nu = d+\alpha+1$, as
this would give the best bounds in the medium range. However, the upper bounds
for $\phi$ established in the previous two sections are too weak to ensure that
$w-\vert v \vert >0$ outside of $\Omega_1(t)$ for all $t$ with that naive
choice. In the long range, we only have $ \vert v(\xi,t) \vert <
e^{\lambda (1-\alpha)t} \left| \xi \right|^{-(d+\alpha)}$. In the short range,
we have that $ \vert v(\xi, t) \vert < e^{(\lambda-\delta_0)t}$ provided the
same can be said of $v$ globally.\\
\\
The added growth rate (in the form of $e^{\eta t}$) alows us to ensure that $w >
\vert v \vert$ in the long range. Since both $w$ and the upper bound for $v$ are
radially decreasing and $O(\left| \xi \right|^{-(d+\alpha)})$, it suffices to
check the inequality for $\left| \xi \right| = ce^{\lambda \alpha t}$. At the
boundary between medium range and long range, $w \approx e^{-\nu \lambda \alpha
t + \eta t}$ and $ \vert v \vert  < e^{-(d+\alpha)\lambda \alpha t + \lambda t -
\lambda \alpha t}$. So it suffices to have $\eta - \lambda \alpha (\nu - d - \alpha) =
\eta - \lambda \alpha (1-\frac{\eta}{\lambda}) \geq \lambda (1 - \alpha)$, or
$\eta \geq \frac{\lambda}{1+\alpha}$. This leaves enough room for $\eta$ to be smaller
than $\lambda$. Fixing constants (and looking at the equation after a transient
period to alow the long range decay estimates to come into effect) we have that
$w > \vert v \vert $ for $\left| \xi \right| > ce^{\lambda \alpha t}$.\\
\\
For the second derivatives, we will similarly employ a supersolution of the form
$W(\xi, t) = e^{\zeta t} \tilde{w}_\nu(\vert \xi \vert, t)$ with nearly the same
$\tilde{w}_\nu$ as in \eqref{W}; specifically, we require a larger constant
$\bar{C}_\nu$ in equation \eqref{odeW}. Essentially, we need $W$ to satisfy
\begin{equation}
\label{WW}
-\partial_t W + \lambda \xi \cdot \nabla W + 2 \lambda W + f'(u) W \leq
-e^{-\lambda \alpha t} \vert \Lambda^\alpha W \vert - Ce^{2 \lambda t-2 \delta t} u^2,
\end{equation}
once again appealing to the full estimate for $\phi$. As before, $\vert
\Lambda^\alpha W(\xi, t) \vert \leq \bar{C}_\nu e^{\zeta t}  \vert \xi
\vert^{-(d+\alpha)}$. But the last term is dominated by $C e^{2 \lambda -2 \delta t} \vert
\xi \vert^{-2(d+\alpha)}$. If we assume that $\zeta > 2 \lambda - 2 \delta$,
the diffusive term is then a stronger perturbation than
the inhomogeneity, and so the latter can be absorbed into the bounds for the
former. Our candidate supersolution $W$ will solve \eqref{WW}
provided that
$$ y \partial_y \tilde{w}_\nu (y,t) + (d+\alpha + 2 - \frac{\zeta}{\lambda})
\tilde{w}_\nu (y,t) \leq \frac{\bar{C}_\nu}{y^{d+\alpha}}. $$
This is just \eqref{odeW} with $\nu = d+\alpha + 2 + \zeta / \lambda$, so the
equation is satisfied. We want to take $\zeta$ as small as possible
while still observing the requirement that $0< \nu-d-\alpha \leq 1$ (this is
necessary to ensure $W$ remains positive on $\Omega_1(t)$ with uniform bounds).
Thus our arguments impose the constraint that $2 \lambda - 2 \delta < \zeta < 2 \lambda$.\\
\\
For the final constraint, $W$ must dominate $V$ (that is, $e^{2 \lambda t} \psi$
in exponential coordinates) at the interface between medium and long ranges. As
before, we have that $W \approx e^{-\nu \lambda \alpha t + \zeta t}$ and $\vert
V \vert \leq e^{-(d+\alpha)\lambda \alpha t + 2 \lambda t - \lambda \alpha t}$ when
$\vert \xi \vert = ce^{\lambda \alpha t}$. We therefore require that $\zeta -
\lambda \alpha (\nu - d - \alpha) = \zeta - \lambda \alpha
(2-\frac{\zeta}{\lambda}) \geq \lambda (2 - \alpha)$, or $\zeta \geq \lambda
\frac{2+\alpha}{1+\alpha}$. For $\alpha > 0$,
this will always overlap with the admissible range of $\zeta$ determined above; and, as
expected, the minimum value for $\zeta$ decreases as we get better diffusion.\\
\\
This establishes, as before, that $W > \vert V \vert$ in the long range and that
$W$ satisfies the comparison inequality on $\Omega_1(t)$. To satisfy all the
conditions for the selective comparison principle, we must show that the
supersolutions for $v$ and $V$ also dominate in the short range, $\Omega_0(t)$.

\subsection{The estimate for $\phi$ and $\psi$}
\noindent{\sc Proof of Theorem \ref{t1.1}.} Recall that $\Omega_0(t) = \lbrace
\vert \xi \vert < \frac{3}{2} \beta \rbrace$ and $\Omega_1(t) = \lbrace \beta <
\vert \xi \vert < c e^{\lambda \alpha t} \rbrace$. Section 3.2  established that, in
$\Omega_0(t)$, $v < C(v_0) e^{(\lambda-\delta_0)t}$ or else $v$ is
larger outside of this regime; which, by Lemma 3.1, in fact means $v$ is larger in
$\Omega_1(t)$; outside of $\Omega_0(t) \cup \Omega_1(t)$ (a compact set), the
comparative exponential bound holds unconditionally.\\
\\
\noindent{}Fix $\eta = \text{max}\left( \frac{\lambda}{1+\alpha}, \lambda -
\delta_0 /2 \right)$ and $w_0$ so large that the conditional bound in the short
range is stronger than the requirement that $w > v$ in the same regime; that is,
$C(v_0) e^{(\lambda-\delta_0)t} < \text{inf}_{\xi \in \Omega_0(t)} w(\xi,t)$
for all $t>0$ (\eqref{W} shows that such a choice is always possible). We can
also ensure that $w > v$ at $t=0$.\\
\\
\noindent{}We then argue by contradiction. Let $\bar{t}>0$ be the critical time such that
$v(\xi,t) < w(\xi,t)$ for all $t< \bar{t}$ and $ \xi \in \Omega_0(t) \setminus \Omega_1(t)$,
and such that the same upper bound fails for a sequence $t_n \rightarrow \bar{t}^{+}$.
This therefore implies a sequence $\lbrace \xi_n \rbrace \in \Omega_0(t_n) \setminus
\Omega_1(t_n)$ such that $v(\xi_n, t_n) \geq w(\xi_n, t_n) \geq C(v_0) e^{(\lambda-\delta_0)t}$.
By the conditional bound of Section 3.2, there is another sequence $\lbrace \xi'_n \rbrace \in
\Omega_1(t_n) \setminus \Omega_0(t_n)$ such that $v(\xi'_n, t_n) \geq
\text{inf}_{\xi \in \Omega_0(t_n) \setminus \Omega_1(t_n)} w(\xi,t_n)$. By compactness and
the continuity of $v$, we extract a limit point $\bar{\xi}$ with $\vert \bar{\xi} \vert \geq
\frac{3}{2} \beta$ and $v(\bar{\xi}, \bar{t}) \geq \text{inf}_{\vert \xi \vert < \beta}
w(\xi, \bar{t})$.\\
\\
\noindent{}However, for $t<\bar{t}$, the hypothesis of Theorem (\ref{SCP}) was by assumption
valid, and therefore $v < w$ up to $\bar{t}$. Moreover, \eqref{W} immediately shows that there
is a continuous, \emph{always positive} $\epsilon(t)$ such that $\text{sup}_{\vert \xi \vert
\geq \frac{3}{2} \beta} w(\xi, t) < \text{inf}_{\vert \xi \vert < \beta} w(\xi, t) -
\epsilon(t)$. The continuity of $v$ allows us to extend the inequality to $\bar{t}$. Thus
$$ v(\bar{\xi}, \bar{t}) \geq \text{inf}_{\vert \xi \vert < \beta}
w(\xi, \bar{t}) \geq \text{sup}_{\vert \xi \vert
\geq \frac{3}{2} \beta} w(\xi, \bar{t}) + \epsilon(\bar{t}) \geq w(\bar{\xi}, \bar{t}) +
\epsilon(\bar{t}) \geq v(\bar{\xi}, \bar{t}) + \epsilon(\bar{t}), $$
which is evidently a contradiction. Therefore, $w > v$ (and analogously $w > -v$) in
$\Omega_0(t) \setminus \Omega_1(t)$.
The hypothesis of Theorem (\ref{SCP}) is always satisfied, and so $w(\xi,t) >
\vert v(\xi,t) \vert$ for all $\xi$ and $t$. Precisely
the same argument also applies to $V$ with $\zeta = \text{max} \left( \lambda
\frac{2+\alpha}{1+\alpha}, 2 \lambda - 2 \delta, 2 \lambda - \delta_0 \right)$.\\
\\
\noindent{}Putting everything together, we have the following three upper bounds (written
in terms of $\phi$). Recall that $u \approx \left| \xi \right|^{-(d+\alpha)}$ in
the medium range.

$$ \vert \phi(x,t) \vert \leq C_2 u(x,t) e^{-\lambda \alpha t} \text{ for } \left| x
\right| > c_1 e^{\lambda (1+\alpha) t} $$
$$ \vert \phi(x,t) \vert \leq C_1 u(x, t) e^{\eta t - \lambda t} \left| \xi
\right|^{-1+\eta/\lambda} \leq C_1 u(x,t) e^{-\delta_1 t} \text{ for } c_0
e^{\lambda t} \leq \left| x \right| \leq c_1 e^{\lambda (1+\alpha) t} $$
$$ \vert \phi(x,t) \vert \leq C_0 e^{-\delta_0 t} \leq C_0 u(x,t) e^{-\delta_0
t} \text{ for } \left| x \right| < c_0 e^{\lambda t} $$
We also get the corresponding three bounds for $\psi$ (recall that $\psi = e^{-2
\lambda t} V$).
$$ \vert \psi(x,t) \vert \leq C_2' u(x,t) e^{-\lambda \alpha t} \text{ for } \left| x
\right| > c_1 e^{\lambda (1+\alpha) t} $$
$$ \vert \psi(x,t) \vert \leq C_1' u(x, t) e^{\zeta t - 2 \lambda t} \left| \xi
\right|^{-2+\zeta/\lambda} \leq C_1' u(x,t) e^{-\bar{\delta_1} t} \text{ for }
c_0 e^{\lambda t} \leq \left| x \right| \leq c_1 e^{\lambda (1+\alpha) t} $$
$$ \vert \psi(x,t) \vert \leq C_0 e^{-\delta_0 t} \leq C_0 u(x,t) e^{-\delta_0
t} \text{ for } \left| x \right| < c_0 e^{\lambda t} $$
Thus proving \eqref{GDE} and Theorem (\ref{t1.1}). \hfill $\bullet$ \\
\\
\noindent{}\textbf{Remark:} Theorem (\ref{t1.1}) can be extended to a comparative exponential
decay for \emph{all} derivatives of $u$. The dual proof presented above is inductive, essentially
repeating the same arguments for $\psi$ while assuming the full result for $\phi$. This can be
reworked into a full induction showing the much stronger result that
$$ | \nabla^k u (x,t) | \leq C u(x,t) e^{-\delta_k t} $$
for all $k$. The next section will, however, only require \eqref{GDE}.
\SE{Bound on the fractional Laplacian}
We finally come to the estimate telling us that, at large times, the fractional
Laplacian of the solution $u$ is small compared to the size of $u$.
\begin{lem}
\label{l4.2}
There exist $C >0$ and $\delta'>0$ such that
\begin{equation}
\vert \Lambda^\alpha u(x,t)\vert  \leq C u(x,t) e^{-\delta' t}.  \label{FDE1}
\end{equation}
\end{lem} 

\noindent{\sc Proof.} We will examine the explicit integral representation for
$\Lambda^\alpha u (x, t)$ and initially with $\left| x \right| \geq c e^{\lambda
t}$. Up to a constant, we have
\begin{equation}
\label{Lambda}
\begin{array}{rll}
\Lambda^\alpha u(x,t) = &\displaystyle\biggl(\int_{\left| x - y \right| <
\epsilon e^{\gamma \lambda t}} dy+ \int_{\left| x - y \right| \geq \epsilon
e^{\gamma \lambda t}} dy\biggl)\frac{u(x,t)-u(y,t)}{\left| x-y
\right|^{d+\alpha}}, \\
= &I_1 + I_2
\end{array}
\end{equation}

\noindent{}where $\gamma \in (0,1]$. For the inner piece $I_1$, we make use of
(\ref{GDE}) and a Taylor expansion: $u(x,t) - u(y,t) = -\nabla u(x,t) \cdot
(y-x) + O \left( \vert \nabla^2 u (cy + (1-c)x,t) \vert \vert x-y \vert^2
\right)$ for some $c \in (0,1)$. Keeping in mind that
$$ \int_A \frac{x-y}{\vert x-y \vert^{d+\alpha}} dy = 0 $$
on any annulus $A$ centered at $x$, we then have
$$ \left| I_1 \right| \leq \int_{\left| x-y \right| < \epsilon e^{\gamma \lambda
t}} \frac{\left| \nabla^2 u (cx+(1-c)y, t) \right|}{\left| x-y
\right|^{d+\alpha-2}} dy \leq e^{-\delta t} \int_{\left| x-y \right| < \epsilon
e^{\gamma \lambda t}} \frac{u(cx+(1-c)y,t)}{\left| x-y \right|^{d+\alpha-2}}  dy.
$$

\noindent{}Now, in this regime, $u(x,t) \approx e^{\kappa t} \left| x
\right|^{-(d+\alpha)}$, and $\left| cx + (1-c)y \right| \geq (1-\epsilon) \left|
x \right|$ because we assumed $\left| x-y \right| < \epsilon e^{\gamma \lambda
t} \leq \epsilon \left| x \right|^\gamma \leq \epsilon \left| x \right|$.
Therefore,
$$ \left| I_1 \right| \leq (1-\epsilon)^{-(d+\alpha)} u(x,t) e^{-\delta t}
\int_0^{\epsilon e^{\gamma \lambda t}} \frac{r^{d-1}dr}{r^{d-\alpha-2}} dr. $$

\noindent{}We then have (up to a constant depending on $\epsilon$)
\begin{equation}
\left| I_1 \right| \leq u(x,t) e^{-\delta t} e^{\gamma \lambda (2-\alpha) t}
\leq u(x,t) e^{-\delta t /2}.
\label{I1}
\end{equation}
To obtain the last inequality we must choose $\gamma$ appropriately. If
$(2-\alpha)\lambda \leq \delta / 2$, we may take $\gamma = 1$. Otherwise, we
take $\gamma = \delta / (2 \lambda (2-\alpha)) < 1$. Either way, we obtain our
final estimate for $I_1$.\\
\\
For the outer piece $I_2$, we merely use the positivity of $u$: $ ( u(x,t) -
u(y,t) ) < u(x,t)$.
$$ \left| I_2 \right| \leq u(x,t) \int_{\left| x-y \right| \geq \epsilon
e^{\gamma \lambda t}} \frac{dy}{\left| x-y \right|^{d+\alpha}} = u(x,t)
\int_{\epsilon e^{\gamma \lambda t}}^\infty \frac{r^{d-1}dr}{r^{d+\alpha}} =
u(x,t) \epsilon^{-\alpha} e^{-\alpha \gamma \lambda t}. $$

\noindent{}The above inequality and (\ref{I1}) show that $\left| \Lambda^\alpha
u(x,t) \right| \leq u(x,t)e^{-\delta' t}$ for some fixed positive $\delta'$ and
$\left| x \right| > ce^{\lambda t}$. The exact value (in terms of $\delta$,
$\alpha$, and $\lambda$) requires an optimization argument for $\gamma$.\\

\noindent{}The case where $\left| x \right| < ce^{\lambda t}$, i.e. the short
range, is much simpler. In this range, we have that $0<\bar{c}<u(x,t) \leq 1$.
And, moreover, we know from Theorem \ref{t1.1} that $\| \nabla^2 u \|_{L^\infty}
\leq Ce^{-\delta t}$. Using our Taylor expansion again, we can interpolate
between these two to obtain:
$$ \left| \int_{{\mathbb R}^d} \frac{u(x,t)-u(y,t)}{\left| x - y
\right|^{d+\alpha}} dy \right| \leq 
\int_{\left| x-y \right| < R} \frac{\| \nabla^2 u (.,t) \|_{L^\infty} dy}{\left|
x-y \right|^{d+\alpha-2}} +
\int_{\left| x-y \right| \geq R} \frac{2 dy}{\left| x-y \right|^{d+\alpha}} \leq
C e^{-\alpha \delta t / 2}. $$
The last inequality required us to optimize $R$ to $\| \nabla^2 u(.,t)
\|_{L^\infty}^{-1/2}$ and use our uniform estimate. Hence $\left|
\Lambda^\alpha u(x,t) \right| \leq Ce^{-\delta' t} \leq C \bar{c}^{-1} u(x,t) e^{-\delta' t}$
for all $x$. \hfill $\bullet$

\SE{Asymptotic symmetrization}
We have used (\ref{GDE}) to establish a similar comparative decay
on the asymptotic size of $\Lambda^\alpha u$. We will now look at the
explicit (integral) representation for the solution $u(x,t)$ and obtain a better
approximate asymptotic behavior than what may be inferred from (\ref{u-bounds}).
We will then use this to prove Theorem \ref{t1.2}.

\subsection{Behavior of $u$ at infinity for $t=1$}
The second ingredient in the symmetrization proof is the fact that $u(x,1)$
decays exactly as $\vert x\vert^{-(d+\alpha)}$ for large $x$. This is detailed
in the following lemma.
\begin{lem}
For all $t>0$, there is $k_t>0$ and $\delta>0$ such that 
\begin{equation}
\label{e4.1}
u(x,t)=\frac{k_t}{\vert x\vert^{(d+\alpha)}}+O(\frac1{\vert
x\vert^{d+\alpha+\delta}}).
\end{equation}
\end{lem}
{\sc Proof.} Once again it is a matter of applying Duhamel's formula to $u$:
$$
\begin{array}{rll}
u(x,t)=&\displaystyle\frac{e^{\kappa
t}}{t^{d/\alpha}}\int_{\RR^d}p_\alpha(\frac{x-y}{t^{1/\alpha}})
u_0(y)dy+\int_0^t\int_{\RR^d}\frac{e^{\kappa(t-s)}}{(t-s)^{d/\alpha}}
p_\alpha(\frac{x-y}{(t-s)^{1/\alpha}})u^2(y,s)
dy\\
=&\displaystyle\frac{e^{\kappa
t}}{t^{d/\alpha}}\int_{\RR^d}p_\alpha(\frac{x-y}{t^{1/\alpha}})
u_0(y)dy+\int_0^te^{\kappa(t-s)}D(s,t,x,y) ds.
\end{array}
$$
We have, by Proposition \ref{p2.1}:
$$
\begin{array}{rll}
\displaystyle \int_{\RR^d}p_\alpha(\frac{x-y}{t^{1/\alpha}})
u_0(y)dy=&\displaystyle c_\alpha t \int_{\RR^d} \frac{u_0(y)}{\vert
x-y\vert^{d+\alpha}}dy+O(\frac1{\vert x\vert^{d+\alpha+\delta}})\\
=&\displaystyle\frac{c_\alpha t}{\vert x\vert^{d+\alpha}} \int_{\RR^d}u_0
dy+O(\frac1{\vert x\vert^{d+\alpha+\delta}})
\end{array}
$$
because $u_0$ is exponentially decreasing. As for the second term, it is
sufficient to estimate $D(s,t,x,y)$, the integration in time being only  between
0 and $t$.
Let us decompose the integral 
$$\displaystyle D(s,t,x,y)=\int_{\vert x-y\vert\leq\vert x\vert/2}+
\int_{\vert x-y\vert\geq\vert x\vert/2}:=I_1+I_2.
$$
We omit, in order to keep the notations light, the dependence in $s$, $t$, $x$,
$y$. As for $I_1$, let us notice the bound 
$
\displaystyle u(x,t)\leq\frac{Ce^{\kappa t}}{\vert x\vert^{d+\alpha}}
$
simply because we have $\partial_t u+\Lambda^\alpha u\leq\kappa u$.
So, we have
$$
\begin{array}{rll}
I_1\leq &\displaystyle\frac{C_t}{\vert
x\vert^{2(d+\alpha)}(t-s)^{d/\alpha}}\int_{\vert x-y\vert\leq\vert
x\vert/2}p_\alpha(\frac{x-y}{(t-s)^{1/\alpha}})dy\\
\leq &\displaystyle\frac{C_t}{\vert x\vert^{2(d+\alpha)}}\Vert
p_\alpha\Vert_{L^1}.
\end{array}
$$
As for $I_2$, let us use estimate \eqref{e2.2} in  Proposition \ref{p2.1} and
write $I_2\leq I_{21}+CI_{22}$, where
$$
I_{22}=(t-s)^{1+\delta/\alpha}\int_{\vert x-y\vert\geq\vert
x\vert/2}\frac{u^2(y,s)}{\vert x-y\vert^{d+\alpha+\delta}}dy\leq\frac{C_t}{\vert
x\vert^{d+\alpha+\delta}},
$$
and
$$
I_{21}=(t-s)\int_{\vert x-y\vert\geq\vert x\vert/2}\frac{u^2(y,s)}{\vert
x-y\vert^{d+\alpha}}dy.
$$
Let us split $I_{21}$ as 
$$
I_{21}=\int_{\vert x-y\vert\geq\vert x\vert/2,\vert y\vert\leq\e\vert
x\vert}+\int_{\vert x-y\vert\geq\vert x\vert/2,\vert y\vert\geq\e\vert
x\vert}=I_{211}+I_{212},
$$
where $\e>0$ is small. We have
$
\displaystyle I_{212}\leq\frac{C_{t,\e}}{\vert x\vert^{2(d+\alpha)}},
$
so it remains to study the last term. We use the fact that, for $\vert
y\vert\leq\e\vert x\vert$
we have
$$
\vert x-y\vert=\vert x\vert(1+O\left(\frac{\vert y\vert}{\vert x\vert}\right) );
$$
in other words we have, setting 
$$
C_x:=\{\vert x-y\vert\geq\vert x\vert/2,\vert y\vert\leq\e\vert x\vert\}
$$
and using once again estimate \eqref{e2.2}:
$$
\begin{array}{rll}
I_{211}=&\displaystyle\frac{c_\alpha(t-s)}{\vert
x\vert^{d+\alpha}}\int_{C_x}u^2(y,s)
dy
+O\left(\frac{t-s}{\vert x\vert^{d+\alpha+1}}\int _{C_x}\vert y\vert u^2(y,s)
dy\right)\\
&+ \displaystyle O\left(\frac{(t-s)^{1+\delta/\alpha}}{\vert
x\vert^{d+\alpha+\delta}}\int_{C_x}u^2(y,s)
dy\right).
\end{array}
$$
Set $l(t)=\displaystyle\int_{\RR^d}u^2(y,t)$ and note that
$\displaystyle\int_{\RR^d}\vert y\vert u^2(y,t)dy<+\infty$.
Hence there is $l_s>0$ such that
$$
I_{211}=\displaystyle\frac{c_\alpha(t-s)l_s}{\vert
x\vert^{d+\alpha}}+O(\frac1{\vert x\vert^{d+\alpha+\delta}}).
$$
Gathering everything yields \eqref{e4.1}.  \hfill $\bullet$
\subsection{Analysis of long-time level sets}
The proof of Theorem \ref{t1.2} now reduces to the explicit resolution of a
simple ODE. 
\noindent{\sc Proof of Theorem \ref{t1.2}.} For all $\lambda\in[0,1)$, there
exists $m_\lambda>0$ such that
$$
\forall u\in[0,\lambda],\  \  m_\lambda u\leq f(u)/\kappa \leq u.
$$
The combination of this remark and Lemma \ref{l4.2} reduces the full problem
\eqref{FKPP} to 
\begin{equation}
\label{e4.5}
\partial_tu=(1+O(e^{-\delta t}))f(u),\   \   \   u(x,1)=\frac{k_1}{\vert
x\vert^{d+\alpha}}+O(\frac1{\vert x\vert^{d+\alpha+\delta}})\ \hbox{as $\vert
x\vert \to+\infty$},
\end{equation}
keeping in mind that $O(e^{-\delta t})$ might depend on $\lambda$. Set 
$$
G(u)=\int\frac{du}{f(u)},
$$
as $u\to 0$ we have
$$
G(u)=\frac{\Log u}{\kappa}+g_0+O(u),
$$
where $g_0$ is a positive constant. Finding the $\lambda$-level set of $u$ for
large $t$ amounts to integrating \eqref{e4.5}
with $x$ large, which yields
$$
\int_1^t  \frac{\partial_t u}{f(u)}=1+O(e^{-\delta t}),
$$
in other words
$$
\Log u+O(u)-\Log\left(\frac{k_1}{\vert x\vert^{d+\alpha}}+O(\frac1{\vert
x\vert^{d+\alpha+\delta}})\right)=\kappa t+q_\infty(u_0)+O(e^{-\delta t}).
$$
Specializing $u(x,t)=\lambda$ yields
$$
{\mathrm{exp}}\left(O(\frac1{\vert
x\vert^{d+\alpha}})\right)\left(\frac{k_1}{\vert
x\vert^{d+\alpha}}+O(\frac1{\vert x\vert^{d+\alpha+\delta}})\right)=
e^{-\kappa t+q_{\lambda,\infty}(u_0)+O(e^{-\delta t})},
$$
for some possibly different constant $q_{\lambda,\infty}(u_0)$. This, after
elementary computations, yields \eqref{e1.5}. \hfill$\bullet$\\
\\
\noindent Let us point out again that we have proved, in quite a strong sense,
that the large time dynamics of the KPP problem  \eqref{FKPP} is the same as
that 
of the ODE $\dot u=f(u)$, the fractional diffusion being there only,  but this
is important to set the initial datum right. This fact was already noticed in
\cite{CCR} and, in
an even more evident fashion, in \cite{MM}.
\section*{Acknowledgment}
This work was initiated during the visit of the second author to the
Mathematical Institute of Toulouse under the NSF GROW program.
The research leading to these results has received funding from the European
Research Council under the European Union's Seventh Framework Program 
(FP/2007-2013) / ERC Grant Agreement n.321186 - ReaDi - Reaction-Diffusion
Equations, Propagation and Modeling. The second author was also partially
supported by the NSF GRFP grant. The authors thank H. Berestycki and X.
Cabré for
raising the problem, A.-C. Coulon for allowing us to reproduce her numerics, and 
A. Zlato$\breve{\text{s}}$ for valuable feedback.
\end{document}